\documentclass[12pt]{amsart}
\usepackage[all]{xy}
\input epsf
\swapnumbers
\newtheorem{thm}[subsection]{Theorem}
\newtheorem{prop}[subsection]{Proposition}

\newenvironment{proo}{\begin{trivlist} \item{\emph{Proof.}}}
  {\hfill $\square$ \end{trivlist}}

\def\DD{\Delta}
\def\aa{\alpha}
\def\bb{\beta}
\def\cc{\gamma}

\def\ss{\sigma}

\def\oo{\omega}
\def\PP{{\mathcal P}}

\def\RRR{{\mathbb R}}

\def\KK{{\mathcal K}}

\def\arbreA{\xymatrix@R=3pt@C=3pt{
&& \\
&*{}\ar@{-}[ur] \ar@{-}[ul] \ar@{-}[d]     &\\
&&
}}

\def\arbreBA{\xymatrix@R=2pt@C=2pt{
&&&&\\
&&&*{}\ar@{-}[ul] & \\
&&*{}\ar@{-}[uurr] \ar@{-}[uull] \ar@{-}[d]     &&\\
&&&&
}}

\def\arbreAB{\xymatrix@R=2pt@C=2pt{
&&&&\\
&*{}\ar@{-}[ur] &&& \\
&&*{}\ar@{-}[uurr] \ar@{-}[uull] \ar@{-}[d]     &&\\
&&&&
}}

\def\arbreBB{\xymatrix@R=2pt@C=2pt{
&&*{}&&\\
&&&& \\
&&*{}\ar@{-}[uurr] \ar@{-}[uull] \ar@{-}[d] \ar@{-}[uu]     &&\\
&&&&
}}

\def\arbreABC{\xymatrix@R=1pt@C=1pt{
&&&&&&\\
&*{}\ar@{-}[ur] &&&&& \\
&&*{}\ar@{-}[uurr] &&&&\\
&&&*{}\ar@{-}[uuurrr] \ar@{-}[uuulll] \ar@{-}[d] &&&\\
&&&&&&
}}

\def\arbreAAC{\xymatrix@R=1pt@C=1pt{
&&&&&&\\
&*{}&&&&& \\
&&*{}\ar@{-}[uurr]\ar@{-}[uu]  &&&&\\
&&&*{}\ar@{-}[uuurrr] \ar@{-}[uuulll] \ar@{-}[d] &&&\\
&&&&&&
}}

\def\arbreBAC{\xymatrix@R=1pt@C=1pt{
&&&&&&\\
&&&*{}\ar@{-}[ul] &&& \\
&&*{}\ar@{-}[uurr] &&&&\\
&&&*{}\ar@{-}[uuurrr] \ar@{-}[uuulll] \ar@{-}[d] &&&\\
&&&&&&
}}

\def\arbreACA{\xymatrix@R=1pt@C=1pt{
&&&&&&\\
&*{}\ar@{-}[ur] &&&&*{}\ar@{-}[ul] & \\
&&&&&&\\
&&&*{}\ar@{-}[uuurrr] \ar@{-}[uuulll] \ar@{-}[d] &&&\\
&&&&&&
}}

\def\arbreACC{\xymatrix@R=1pt@C=1pt{
&&&&&&\\
&*{}\ar@{-}[ur]&&&&&\\
&&&&& & \\
&&&*{}\ar@{-}[uuurrr] \ar@{-}[uuulll] \ar@{-}[d] \ar@{-}[uuu]&&&\\
&&&&&&
}}

\def\arbreCAB{\xymatrix@R=1pt@C=1pt{
&&&&&&\\
&&&*{}\ar@{-}[ur] &&& \\
&&&&*{}\ar@{-}[uull] &&\\
&&&*{}\ar@{-}[uuurrr] \ar@{-}[uuulll] \ar@{-}[d] &&&\\
&&&&&&
}}

\def\arbreCBA{\xymatrix@R=1pt@C=1pt{
&&&&&&\\
&&&&&*{}\ar@{-}[ul] & \\
&&&&*{}\ar@{-}[uull] &&\\
&&&*{}\ar@{-}[uuurrr] \ar@{-}[uuulll] \ar@{-}[d] &&&\\
&&&&&&
}}

\def\lacet{\xymatrix{
&&*{}&&\\
&&&&\\
&*{}\ar@/^2pc/@{-}[ruu]&&*{}\ar@/_2pc/@{-}[luu]&\\
a=&&&&\\
&&*{}\ar[luu] \ar@{-}[ruu]&\!\! x_0&
}}

\def\deuxlacets{\xymatrix{
&&&&& &&&&& \\
&*{}&&&& &&&&*{}& \\
&&&&*{}\ar@/_2pc/@{-}[lllu]& &*{}\ar@/^2pc/@{-}[rrru]&&&& \\
&a&&&&&&&&b&\\
&&*{}\ar@/^4pc/@{-}[luuu]&&&
*{}\ar[lll]\ar@{-}[luu]\ar[ruu]\ar@{-}[rrr]
 &&&*{}\ar@/_4pc/@{-}[ruuu]&& \\
&&&&&x_0 &&&&&\\
}}

\def\Kzero{\xymatrix@R=4pt@C=4pt{
\\
\\
\\
{\bullet}\\
}}

\def\KzeroC{\xymatrix@R=4pt@C=4pt{
\\
\\
(1)\\
{\bullet}\\
}}

\def\KunA{\xymatrix@R=4pt@C=4pt{
&&\\
&&\\
&&\\
*{}\ar@{-}[rr]&&*{}\\
}}

\def\KunAprim{\xymatrix@R=4pt@C=4pt{
&&*{}\\
*{}\ar@{-}[urr]&&\\
}}

\def\Kungrossi{\xymatrix@R=4pt@C=4pt{
&&&&&\\
&&*{}\ar@{-}[dll]\ar@{-}[ddrrr]  &&& \\
*{}\ar@{-}[dd] &&&&&\\
&&&&&*{} \\
*{}&&&&& \\
&&*{}&&& \\
}}

\def\Kungrossicone{\xymatrix@R=4pt@C=4pt{
&&&&&\\
&&*{}\ar@{-}[dll]\ar@{-}[ddrrr]  &&& \\
*{}\ar@{-}[dd] &&&&&\\
&&&&&*{}\ar@{.}[ddlll]  \\
*{}\ar@{.}[drr] &&&&& \\
&&*{}&&&\\
}}

\def\KunC{\xymatrix@R=4pt@C=4pt{
&&\\
&&\\
(1,2)&&(2,1)\\
*{}\ar@{-}[rr]&&*{}\\
}}

\def\KunF{\xymatrix@R=4pt@C=4pt{
&&\\
&&\\
&&\\
*{}\ar@{->}[rrr]&&&*{}\\
}}

\def\KunSquare{\xymatrix@R=20pt@C=20pt{
&&&&\\
*{}\ar@{|-}[r]&*{}\ar@{|-|}[r]&&,&\\
}}

\def\TdeuxA{\xymatrix@R=4pt@C=4pt{
&&*{}\ar@{-}[dddll]\ar@{-}[dddrr]&&\\
&&&&\\
&&&&\\
*{}\ar@{-}[rrrr]&&&&*{}\\
}}

\def\GdeuxA{\xymatrix@R=4pt@C=4pt{
&&&&\\
&&&&\\
*{}\ar@{-}@/^1pc/[rrrr]\ar@{-}@/_1pc/[rrrr]&&&&*{}\\
&&&&\\
&&&&\\
}}

\def\TtroisA{\xymatrix@R=4pt@C=4pt{
&&&*{}\ar@{-}[ddddlll]\ar@{-}[ddddrrr]\ar@{-}[drrr]&&&\\
&&&&&&*{}\ar@{-}[ddd]\ar@{.}[dddllllll]\\
&&&&&*{}&\\
&&&&&&\\
*{}\ar@{-}[rrrrrr]&&&&&&*{}\\
}}

\def\GtroisA{\xymatrix@R=4pt@C=4pt{
&&&*{}\ar@{-}@/^0.5pc/[ddrrr]&&&\\
&&&&&&\\
*{}\ar@{.}@/^1pc/[rrrrrr]\ar@{-}@/_1pc/[rrrrrr]\ar@{-}@/^0.5pc/[uurrr]\ar@{-}@/_0.5pc/[ddrrr]&&&&&&*{}\\
&&&&&&\\
&&&*{}\ar@{-}@/_0.5pc/[uurrr]&&&\\
}}

\def\TquatreA{\xymatrix@R=4pt@C=4pt{
&&&*{}\ar@{-}[ddddlll]\ar@{-}[ddddrrr]\ar@{-}[drrr]&&&\\
&&&&&&*{}\ar@{-}[ddd]\ar@{.}[dddllllll]\\
&&&&&*{}&\\
&&&&*{}\ar@{.}[dllll]\ar@{.}[uuul]\ar@{.}[uurr]\ar@{.}[drr]&&\\
*{}\ar@{-}[rrrrrr]&&&&&&*{}\\
}}

\def\QdeuxA{\xymatrix@R=3pt@C=3pt{
*{}\ar@{-}[rrrr]\ar@{-}[dddd]&&&&*{}\ar@{-}[dddd]\\
&&&&\\
&&&&\\
&&&&\\
*{}\ar@{-}[rrrr]&&&&*{}\\
}}

\def\QtroisA{\xymatrix@R=5pt@C=5pt{
&&*{}\ar@{-}[rrr]\ar@{.}[ddd]&&&*{}\ar@{-}[ddd]\\
*{}\ar@{-}[rrr]\ar@{-}[ddd]\ar@{-}[urr]&&&*{}\ar@{-}[ddd]\ar@{-}[urr]&&\\
&&&&&\\
&&*{}\ar@{.}[rrr]&&&*{}\\
*{}\ar@{-}[rrr]\ar@{.}[urr]&&&*{}\ar@{-}[urr]&&\\
}}

\def\KdeuxA{\xymatrix@R=4pt@C=4pt{
&&&&&\\
&&*{}\ar@{-}[dll]\ar@{-}[ddrrr]  &&& \\
*{}\ar@{-}[dd] &&&&&\\
&&&&&*{}\ar@{-}[ddlll]  \\
*{}\ar@{-}[drr] &&&&& \\
&&*{}&&& \\
}}

\def\KdeuxAprim{\xymatrix@R=0pt@C=0pt{
& &&&&&& && \\
&*{}\ar@{-}[rrrrr]\ar@{-}[ldd]  &&& &&*{}\ar@{-}[ddddll]& && \\
& &&& &&& && \\
*{}\ar@{-}[rdd] & &&& &&& & &\\
& &&&&&& && \\
&*{}\ar@{-}[rrr] &&&*{} &&& & &\\
}}

\def\KdeuxAgrossi{\xymatrix@R=0pt@C=0pt{
&*{}\ar@{-}[rrrrrrrrrr] *{}\ar@{-}[rrrrddd] *{}\ar@{-}[ldd] &&&& &&& &&&*{}\ar@{-}[rrd] && \\
&&&& &&& &&& &&&*{}\ar@{-}[llldd] *{}\ar@{-}[llldddddd]   \\
*{}\ar@{-}[rrrrddd]  *{}\ar@{-}[rrdddd] &&&&& &&& &&& && \\
&&&&&*{}\ar@{-}[rrrrr] *{}\ar@{-}[ldd]  &&& &&*{}\ar@{-}[lldddd] & && \\
&&&&& &&& &&& && \\
&&&&*{}\ar@{-}[rdd] & &&& &&& & &\\
&&*{}*{}\ar@{-}[rrrddd] &&& &&&*{} &&& && \\
&&&&&*{}\ar@{-}[rrr] *{}\ar@{-}[dd]  &&&*{}\ar@{-}[dd]  &&*{}\ar@{-}[lldd] & & &\\
&&&&& &&& &&& && \\
&&&&&*{}\ar@{-}[rrr]  &&&*{} &&& && \\
}}

\def\KdeuxSquare{\xymatrix@R=4pt@C=4pt{
& &*{}\ar@{-}[dl] \ar@{-}[dr]& &  &\\
&*{}\ar@{-}[dl]\ar@{-}[ddr] & &*{}\ar@{-}[ddl]\ar@{-}[ddrr]] &  &\\
*{}\ar@{-}[d]& & & &  &\\
*{}\ar@{-}[d] \ar@{-}[rr]& &*{}\ar@{-}[ddl] \ar@{-}[ddr]& &  &*{}\ar@{-}[ddll]\\
*{}\ar@{-}[dr]& & & &  &\\
& *{}\ar@{-}[dr]& & *{}\ar@{-}[dl] &  &\\
& & *{}& &  &\\
}}
\def\KdeuxC{\xymatrix@R=4pt@C=4pt{
&&&(1,2,3)&&&\\
&&&*{}\ar@{-}[dll]\ar@{-}[ddrrr]  &&&& \\
(2,1,3)&*{}\ar@{-}[dd] &&&&&&\\
&&&&&&*{}\ar@{-}[ddlll]&(1,4,1)  \\
(3,1,2)&*{}\ar@{-}[drr] &&&&&& \\
&&&*{}&&&& \\
&&&(3,2,1)&&&&
}}

\def\KdeuxF{\xymatrix@R=4pt@C=4pt{
&& &&\\
&&*{}\ar[dll]\ar[ddrrr]  &&&&\\
 *{}\ar[dd] &&&&&\\
&&&&&*{}\ar[ddlll] \\
 *{}\ar[drr] &&&&& \\
&&*{}&&& \\
&& &&&
}}

\def\KdeuxT{\xymatrix@R=10pt@C=10pt{
&&*{}\ar[dll]\ar[ddrrr]\ar[dddd] &&&&\\
 *{}\ar[dd] \ar[dddrr]& &&&&\\
&&& &&*{}\ar[ddlll] \\
 *{}\ar[drr] &&&&& \\
&&*{}&&& \\
}}

\def\KdeuxTbis{\xymatrix@R=10pt@C=10pt{
&&*{} \ar[dddll]\ar[dll]\ar[ddrrr]\ar[dddd] &&&&\\
 *{}\ar[dd]& &&&&\\
&&& &&*{}\ar[ddlll] \\
 *{}\ar[drr] &&&&& \\
&&*{}&&& \\
}}

\def\KdeuxFindic{\xymatrix@R=10pt@C=10pt{
&&*{}\ar@{=>}[dll]\ar@{=>}[ddrrr] \ar[dddd] &&&&\\
 *{}\ar@{=>}[dd] \ar[dddrr]& 12 &&&&\\
&&& 21 &&*{}\ar[ddlll] \\
 *{}\ar[drr] &11&&&& \\
&&*{}&&& \\
}}

\def\KdeuxFindicbis{\xymatrix@R=10pt@C=10pt{
&&{}^{a=1}_{p=0}&*{}\ar@{=>}[dll]\ar@{=>}[ddrrr] \ar[dddd] &&&&\\
& *{}\ar@{=>}[dd] \ar[dddrr]& 12 &&&{}^{a=2}_{p=1}&\\
{}^{a=1}_{p=1}&&&& 21 &&*{}\ar[ddlll] \\
&*{}\ar[drr] &11&&&& \\
&&&*{}&&& \\
}}

\def\KdeuxFindicSimple{\xymatrix@R=10pt@C=10pt{
&&*{}\ar[dll]\ar[ddrrr] \ar[dddd] &&&&\\
 *{}\ar[dd] \ar[dddrr]& 12 &&&&\\
&&& 21 &&*{}\ar[ddlll] \\
 *{}\ar[drr] &11&&&& \\
&&*{}&&& \\
}}

\def\KtroisFtriang{\xymatrix@R=20pt@C=20pt{
&*{}\ar@{=>}[rrrrrrrrrr] *{}\ar@{<=}[rrrrddd] *{}\ar@{=>}[ldd] &&&& &&& &&&*{}\ar@{<=}[rrd] && \\
&&&& &&& &&& &&&*{}\ar@{<=}[llldd] *{}\ar@{=>}[lldddddd]   \\
*{}\ar@{<=}[rrrrddd]  *{}\ar@{=>}[rrdddd] &&&&& &&& &&& && \\
&&&&&*{}\ar@{=>}[rrrrr] \ar@{=>}[ldd] \ar[ulllll] \ar[uuurrrrrr] \ar[ddddrrr] &&& &&*{}\ar@{=>}[lldddd] \ar[uuur] \ar[rdddd]& && \\
&&&&& &&& &&& && \\
&&&&*{}\ar@{=>}[rdd] \ar[dll] \ar[ddrrrr]& &&& &&& & &\\
&&*{}\ar@{<=}[rrrddd] &&& &&&*{}&&& && \\
&&&&&*{}\ar@{=>}[rrr] *{}\ar@{=>}[dd] \ar[ulll] \ar[ddrrr]&&&*{}\ar@{=>}[dd] \ar[rrr] &&&*{}\ar@{<=}[llldd]  & &\\
&&&&& &&& &&& && \\
&&&&&*{}\ar@{=>}[rrr]  &&&*{} &&& && \\
}}

\def\KtroisFindic{\xymatrix@R=12pt@C=12pt{
&*{}\ar@{=>}[rrrrrrrrrr] *{}\ar@{<=}[rrrrddd] *{}\ar@{=>}[ldd] &&&& &&& &&&*{}\ar@{<=}[rrd] && \\
&&&& & 321 && &&& & 311 &&*{}\ar@{<=}[llldd] *{}\ar@{=>}[lldddddd]   \\
*{}\ar@{<=}[rrrrddd]  *{}\ar@{=>}[rrdddd] &&231&&& &&& 312 &&& && \\
&&&&&*{}\ar@{=>}[rrrrr] \ar@{=>}[ldd] \ar[ulllll] \ar[uuurrrrrr] \ar[ddddrrr] &&& &&*{}\ar@{=>}[lldddd] \ar[uuur] \ar[rdddd]& && \\
&&&& 213 & && 132  &&&& 121 && \\
&& 221 &&*{}\ar@{=>}[rdd] \ar[dll] \ar[ddrrrr]& 123 &&& &&& & &\\
&&*{}\ar@{<=}[rrrdddd] && 212 & 122 &&&*{}& 112& & && \\
&&&& 211 &*{}\ar@{=>}[rrr] *{}\ar@{=>}[ddd] \ar[ulll] \ar[dddrrr]&&&*{}\ar@{=>}[ddd] \ar[rrr] &&&*{}\ar@{<=}[lllddd]  & &\\
&&&&& &   & 113 &&  111 && && \\
&&&&&&  131 &  &&    && && \\
&&&&&*{}\ar@{=>}[rrr]  &&&*{} &&& && \\
}}

\def\PdeuxC{\xymatrix@R=2pt@C=2pt{
&&&(1,2,3)&&\\
&&&*{}\ar@{-}[ddll]\ar@{-}[ddrr]  &&& \\
&&&&&&\\
(2,1,3)&*{}\ar@{-}[dd] &&&&*{}\ar@{-}[dd] &(1,3,2)\\
&&&&&  \\
(3,1,2))&*{}\ar@{-}[ddrr] &&&&*{}\ar@{-}[ddll]&(2,3,1) \\
&&&&&&\\
&&&*{}&&& \\
&&&(3,2,1)&&&
}}

\def\PdeuxA{\xymatrix@R=4pt@C=4pt{
&&&&&\\
&&*{}\ar@{-}[dll]\ar@{-}[drr]  &&& \\
*{}\ar@{-}[dd] &&&&*{}\ar@{-}[dd] &\\
&&&&& \\
*{}\ar@{-}[drr] &&&&*{}\ar@{-}[dll] & \\
&&*{}&&& \\
}}

\def\KtroisA{\xymatrix@R=0pt@C=0pt{
&*{}\ar@{-}[rrrrrrrrrr] *{}\ar@{-}[rrrrddd] *{}\ar@{-}[ldd] &&&& &&& &&&*{}\ar@{.}[llldddddd]  *{}\ar@{-}[rrd] && \\
&&&& &&& &&& &&&*{}\ar@{-}[llldd] *{}\ar@{-}[llldddddd]   \\
*{}\ar@{-}[rrrrddd]  *{}\ar@{-}[rrdddd] &&&&& &&& &&& && \\
&&&&&*{}\ar@{-}[rrrrr] *{}\ar@{-}[ldd]  &&& &&*{}\ar@{-}[lldddd] & && \\
&&&&& &&& &&& && \\
&&&&*{}\ar@{-}[rdd] & &&& &&& & &\\
&&*{}\ar@{.}[rrrrrr] *{}\ar@{-}[rrrddd] &&& &&&*{}\ar@{.}[rrd]  &&& && \\
&&&&&*{}\ar@{-}[rrr] *{}\ar@{-}[dd]  &&&*{}\ar@{-}[dd]  &&*{}\ar@{-}[lldd] & & &\\
&&&&& &&& &&& && \\
&&&&&*{}\ar@{-}[rrr]  &&&*{} &&& && \\
}}

\def\KtroisF{\xymatrix@R=0pt@C=0pt{
&*{}\ar[rrrrrrrrrr] *{}\ar@{<-}[rrrrddd] *{}\ar[ldd] &&&& &&& &&&*{}\ar@{.>}[llldddddd]  *{}\ar@{<-}[rrd] && \\
&&&& &&& &&& &&&*{}\ar@{<-}[llldd] *{}\ar[llldddddd]   \\
*{}\ar@{<-}[rrrrddd]  *{}\ar[rrdddd] &&&&& &&& &&& && \\
&&&&&*{}\ar[rrrrr] *{}\ar[ldd]  &&& &&*{}\ar[lldddd] & && \\
&&&&& &&& &&& && \\
&&&&*{}\ar[rdd] & &&& &&& & &\\
&&*{}\ar@{.>}[rrrrrr] *{}\ar@{<-}[rrrddd] &&& &&&*{}\ar@{<.}[rrd]  &&& && \\
&&&&&*{}\ar[rrr] *{}\ar[dd]  &&&*{}\ar[dd]  &&*{}\ar@{<-}[lldd] & & &\\
&&&&& &&& &&& && \\
&&&&&*{}\ar[rrr]  &&&*{} &&& && \\
}}

\def\KtroisFOld{\xymatrix@R=3pt@C=3pt{
&*{}\ar[rrrrrrrrrr] *{}\ar@{<-}[rrrrddd] *{}\ar[ldd] &&&& &&& &&&*{}\ar@{.>}[llldddddd]  *{}\ar@{<-}[rd] & \\
&&&& &&& &&& &&*{}\ar@{<-}[lldd] *{}\ar[llldddddd]   \\
*{}\ar[rrrrddd]  *{}\ar@{-}[rrdddd] &&&&& &&& &&& & \\
&&&&&*{}\ar[rrrrr] *{}\ar[ldd]  &&& &&*{}\ar@{<-}[lldddd] & & \\
&&&&& &&& &&& & \\
&&&&*{}\ar[rdd] & &&& &&& & \\
&&*{}\ar@{.>}[rrrrrr] *{}\ar@{<-}[rrrddd] &&& &&&*{}\ar@{<.}[rd]  &&& & \\
&&&&&*{}\ar[rrr] *{}\ar[dd]  &&&*{}\ar[dd]  &*{}\ar@{<-}[ldd] && & \\
&&&&& &&& &&& & \\
&&&&&*{}\ar[rrr]  &&&*{} &&& & \\
}}

\def\PtroisA{\xymatrix@R=4pt@C=4pt{
&&&*{}\ar@{-}[dll]\ar@{.}[dr]\ar@{-}[r] &*{}\ar@{-}[drr]\ar@{-}[dll] &&& \\
&*{}\ar@{-}[ddl]\ar@{-}[r]&*{}\ar@{-}[dd] & &*{}\ar@{.}[ddl]\ar@{.}[drr]&&*{}\ar@{-}[dd]\ar@{-}[dr]& \\
&&& & &&*{}\ar@{.}[r]\ar@{.}[ddl]&*{}\ar@{-}[dd] \\
*{}\ar@{-}[dd]\ar@{.}[dr]&&*{}\ar@{-}[ddl]\ar@{-}[drr]&*{}\ar@{.}[dll]\ar@{.}[drr] & &&*{}\ar@{-}[dll]\ar@{-}[dr]& \\
&*{}\ar@{.}[dd]&& &*{}\ar@{-}[ddl] &*{}\ar@{.}[dd]&&*{}\ar@{-}[ddl] \\
*{}\ar@{-}[dr]\ar@{-}[r]&*{}\ar@{-}[drr]&& & &&& \\
&*{}\ar@{-}[drr]&&*{}\ar@{-}[dr] & &*{}\ar@{.}[dll]\ar@{.}[r]&*{}\ar@{-}[dll]& \\
&&&*{}\ar@{-}[r] &*{} &&& \\
}}

\def\gammaabctree{\xymatrix@R=3pt@C=3pt{
&&*{}&*{}& &*{}&*{}&& & q+2 \textrm{ leaves}\\
*{}&*{}&&&*{}\ar@{-}[ull]\ar@{-}[ul] \ar@{-}[ur]\ar@{-}[urr]&&*{}& 1 &*{} &\\
p+1&&&& \quad a &&&& 0& \\
&&&& *{}\ar@{-}[uullll]\ar@{-}[uulll] \ar@{-}[uurr]\ar@{-}[uurrrr] \ar@{-}[uu] &&&&& \\
&&&&*{}\ar@{-}[u] &&&&&
}}

\def\squareT{\xymatrix@R=3pt@C=3pt{
(0,1)&*{}\ar[rrrr]&&& &*{}&(1,1) \\
&&&& && \\
&&21&& && \\
&&&&12 && \\
&&&& && \\
(0,0)&*{}\ar[uuuuu]\ar[uuuuurrrr]\ar[rrrr]&&& &*{}\ar[uuuuu]&(1,0)
}}

\def\Pdeuxoriente{\xymatrix@R=4pt@C=4pt{
&&{\quad }&*{}\ar[ddll]\ar[ddrr]\ar[ddddddd]   &{\quad  }&& \\
&&&&&&\\
&*{}\ar[ddd] \ar[dddddrr] &&&&*{}\ar[ddd] \ar[dddddll] &\\
&&&&&&& \\
&&&&&&& \\
&*{}\ar[ddrr] &&&&*{}\ar[ddll] & \\
&&&&&&\\
&&&*{}&&& \\
}}

\begin{document}

\author[J.-L. Loday]{Jean-Louis Loday}
\address{Institut de Recherche Math\'ematique Avanc\'ee\\
    CNRS et Universit\'e Louis Pasteur\\
    7 rue R. Descartes\\
    67084 Strasbourg Cedex, France}
\email{loday@math.u-strasbg.fr}

\title{Parking functions and triangulation of the associahedron}


\begin{abstract}We show that a minimal triangulation of the associahedron (Stasheff polytope) of dimension $n$ is made of $(n+1)^{n-1}$ simplices. We construct a natural bijection with the set of parking functions from a new interpretation of parking functions in terms of shuffles.
\end{abstract}

\maketitle

\section*{Introduction} \label{S:int}
The Stasheff polytope, also known as the associahedron, is a polytope which comes naturally with a poset structure on the set of vertices (Tamari poset), hence a natural orientation on each edge. We decompose this polytope into a union of oriented simplices, the orientation being compatible with the poset structure. This construction defines the associahedron as the geometric realization of a simplicial set. In dimension $n$ the number of (non-degenerate) simplices is $(n+1)^{n-1}$.

A parking function is a permutation of a sequence of integers $i_{1}\leq \cdots \leq i_{n}$ such that $1\leq i_{k}\leq k$ for any $k$. For fixed $n$ the number of parking functions is $(n+1)^{n-1}$. We show that the set $PF_{n}$ of parking functions of length $n$ admits the following inductive description:
$$PF_{n} = \bigcup_{{}^{p+q=n-1}_{p\geq 0, q\geq 0}}\{1,\ldots , p+1\}\times Sh(p,q)\times PF_{p}\times PF_{q}$$
where $Sh(p,q)$ is the set of $(p,q)$-shuffles. From this bijection we deduce a natural bijection between the top dimensional simplices of the associahedron and the parking functions.

$$\KdeuxFindicSimple$$

In the last section we investigate a similar triangulation of the permutohedron.
\medskip

Thanks to Andr\'e Joyal for mentioning to me Abel's formula during the Street's fest.

 \section{Associahedron}\label{S:associahedron}

\subsection{Planar trees, Stasheff complex}\label{Stasheff} The associahedron can be constructed as a cellular complex as follows (cf. for instance \cite{BV}).

Let $Y_n$ be the set of planar binary rooted trees with $n$ internal vertices:
\begin{displaymath}
Y_0 = \{\  \vert \ \}\  , \quad Y_1=\big\{ \vcenter{\arbreA}\big\}\  , \quad 
Y_2=\Big\{\vcenter{\arbreBA}\ ,\  \vcenter{\arbreAB }\Big\}\ ,
\end{displaymath}
\begin{displaymath}
{Y_3=\bigg\{ \vcenter{\arbreABC}}\ ,\  {{}\vcenter{\arbreBAC }},\  {{}\vcenter{\arbreACA }},\  {{}\vcenter{\arbreCAB }},\  {{}\vcenter{\arbreCBA }}\bigg\}\ .
\end{displaymath}
Observe that $t\in Y_{n+1}$ has $n$ internal edges.
For each $t\in Y_{n+1}$ we take a copy of the cube $I^n$ (where $I=[0,1]$ is the interval) which we denote by $I^n_{t}$. Then the associahedron of dimension $n$ is the quotient
$$\KK^n:= \bigsqcup_{t}I^n_{t}/ \sim$$
where the equivalence relation is as follows. We think of an element $\tau=(t; \lambda_{1},\ldots ,\lambda_{n})\in I^n_{t}$ as a tree of type $t$ where the $\lambda_{i}$'s are the lengths of the internal edges. If some of the $\lambda_{i}$'s are 0, then the geometric tree determined by $\tau$ is not binary anymore (since some of its internal edges have been shrinked to a point). We denote the new tree by $\bar{\tau}$. For instance, if none of the $\lambda_{i}$'s is zero, then $\bar{\tau}=t$ ; if all the $\lambda_i$'s are zero, then the tree $\bar{\tau}$ is the corolla (only one vertex). The equivalence relation $\tau \sim \tau'$ is defined by the following two conditions:

-  $\bar{\tau}=\bar{\tau'}$,

-  the lengths of the nonzero-length edges of $\tau$ are the same as those of $\tau'$.

 Hence $\KK^n$ is obtained as a cubical realization:
$$ \quad \KunSquare   \quad \KdeuxSquare
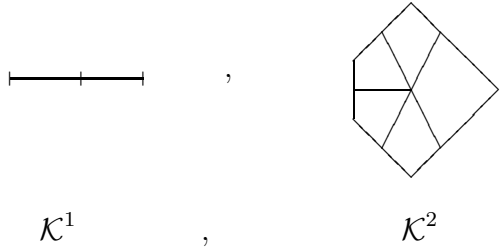
$$
$$\KK^1\qquad \qquad ,\qquad \qquad \qquad   \KK^2 $$

Since a cube can be decomposed into simplices, we can get a simplicial decomposition of $\KK^n$. However our aim is to construct a minimal simplicialization.
It was shown by Stasheff in \cite{Sta} that $\KK^n$ is homeomorphic to a ball. In fact this Stasheff complex can be realized as a polytope (cf. \cite{Lee},\cite{GKZ},\cite{SS},\cite{Lod2}). One way to construct it is recalled in the following section (taken from \cite{Lod2}).

\subsection{Associahedron} \label{tree} To any tree $t\in Y_n, n>0,$ we associate a point $M(t)\in \RRR^n$ with integral coordinates as follows. Let us number the leaves of $t$ from left to right by $0, 1, \ldots , n$. So one can number the vertices from $1$ to $n$ (the vertex  number $i$ is in between the leaves $i-1$ and $i$). We consider the subtree generated by the $i$th vertex. Let $a_i$  be the number of offspring leaves on the left side of the vertex  $i$ and let $b_i$ be the number of offspring leaves on the right side. Observe that these numbers depend only on the subtree determined by  the vertex  $i$. We define:
$$M(t):= (a_1b_1, \ldots , a_ib_i, \ldots ,a_nb_n)\in \RRR^n \ .$$
In low dimension we get:
\begin{eqnarray*}
M(\ \vcenter{\arbreA} )&=& (1\times 1) = (1), \\
 M\big(\ \vcenter{\arbreAB}\big)&=& (1\times 1, 2\times 1) = (1,2), \\
M\big(\ \vcenter{\arbreBA}\big)&=& (1\times 2, 1\times 1)= (2,1),\\
 M\Big(\ \vcenter{\arbreABC}\Big)&=& (1\times 1, 2\times 1, 3\times 1)= (1,2,3), \\
M\Big(\ \vcenter{\arbreACA}\Big) &=&(1\times 1, 2\times 2, 1\times 1)= (1,4,1).\\
\end{eqnarray*}
The planar binary trees with $n$ internal vertices are in bijection with the parenthesizings of the word $x_0x_1\cdots x_{n+1}$. 
For the tree corresponding to the parenthesizing $(((x_0x_1)x_2)(x_3x_4))$  one gets the following point
$(1\times 1,2\times 1, 3\times 2,1\times 1)=(1,2,6,1)$.

Denote by  $H_n$ the hyperplane of $\RRR^n$ whose equation is $$x_1+\cdots +x_n = \frac{n(n+1)}{2} .$$
One can show that  for any tree $t\in Y_n$ the point $M(t)$ belongs to the hyperplane $H_n$.

Then by \cite{Lod2} the \emph{associahedron} or \emph{ Stasheff polytope} $\KK^{n-1}$ is the convex hull of the points $M(t)$ in the hyperplane $H_n$ for $t\in Y_n$.
\\
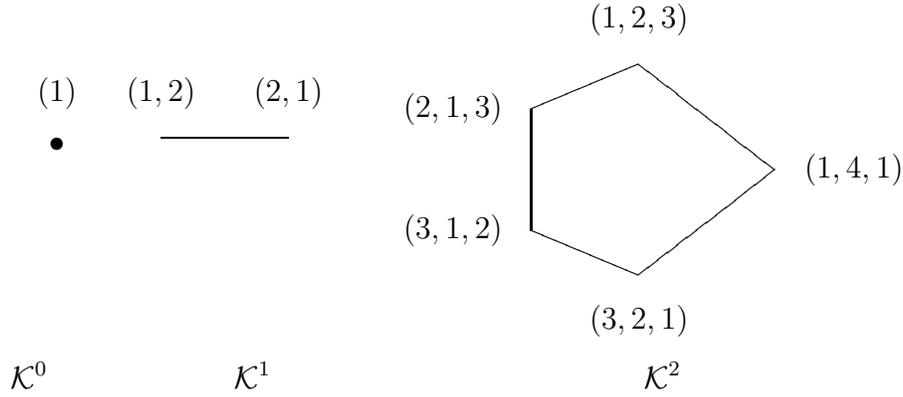
$$\KzeroC\quad  \KunC\qquad \KdeuxC$$
$$\KK^0\qquad\qquad\qquad \KK^1\qquad\qquad\qquad\qquad\qquad\qquad\KK^2\qquad\qquad\qquad\qquad$$

\subsection{Order structure}
Let us recall that, on the set $Y_n$, there is a partial order known as the Tamari order.
It is induced by the order $\vcenter{\arbreAB} \to \vcenter{\arbreBA}$ on $Y_2$ as follows. There is a covering relation $t\to t'$ between two elements of $Y_n$ if $t'$ can be obtained from $t$ by replacing locally a subtree of the form $\vcenter{\arbreAB} $ by a subtree of the form $\vcenter{\arbreBA}$. In low dimension the covering relations induce the following order on the vertices:

$$\KunF\qquad\quad \KdeuxF\qquad\quad \KtroisF$$

Our aim is to triangulate (we mean simplicialize) the associahedron $\KK^n$ by (oriented) $n$-simplices, so that the oriented edges of the simplices are coherent with the Tamari order. We observe immediately that there are two choices for $\KK^2$:

$$\KdeuxT \qquad \KdeuxTbis$$

We choose the first one and we will construct a triangulation for $\KK^n$ consistent with this choice.

\subsection{Boundary of $\KK^n$}\label{boundary} The cells of the associahedron are in bijection with the planar rooted trees (see for instance \cite{LR2}). The vertices correspond to the binary trees and the big cell corresponds to the corolla. The boundary of $\KK^n$, denoted $\partial \KK^n$, is made of cells of the form $\KK^p\times \KK^q$, $p+q=n-1$. They are in bijection with the trees $\cc(a;p,q)$ with two vertices:

$$\gammaabctree$$

Here $p+2$ is the number of outgoing edges at the root, $q+2$ is the number of outgoing edges (leaves) at the other vertex, and $a$ is the index of the only edge which is not a leaf, so $0\leq a \leq p+1$. By convention we index the edges at the root from 0 to $p+1$ from right to left.

Let us denote by $S$ (like South pole) the vertex of $\KK^n$ with coordinates $(n, n-1, \ldots , 1)$, which corresponds to the right comb.

\begin{prop}\label{cone} The associahedron $\KK^n$ is the cone with vertex $S$ over the union of the cells $\cc(a;p,q), a\geq 1$, in $\partial \KK^n$. 
\end{prop}
\begin{proo} 
Since $\KK^n$ is a ball, $\partial \KK^n$ is a sphere. The $n$-cells of $\partial \KK^n$ which contain $S$ are such that $a=0$, because the tree of such a cell is obtained from the right comb (by collapsing $n-2$ edges). The other ones, for which $a\geq 1$, form a ball of dimension $n-1$ and obviously $\KK^n$ can be viewed as a cone with vertex $S$ over this $(n-1)$-dimensional ball.
\end{proo}

\begin{thm}\label{induction} The associahedron $\KK^n$ can be constructed out of $\KK^{n-1}$ as follows:

(a) start with $\KK^{n-1}$,

(b) ``fatten" $\KK^{n-1}$ by replacing its boundary faces $\cc(a-1;p-1,q)$ (of the form $\KK^{p-1}\times \KK^q$), $p+q=n-1$, by  $\cc(a;p,q)$ (of the form $\KK^{p}\times \KK^q$),

(c) take the cone over the fat-$\KK^{n-1}$.
\end{thm}
\begin{proo} From Proposition \ref{cone} it suffices to check that, in $\KK^{n}$, the union of the faces $\cc(a;p,q)$, $a\geq 1$, is precisely fat-$\KK^{n-1}$. Indeed, $\cc(1;0,n-1)$ corresponds to $\KK^{n-1}$ and all the other faces $\cc(a;p,q)$, $a\geq 1, p\geq 1$, come from the cells  $\cc(a-1;p-1,q)$ of $\partial \KK^{n-1}$ by fattening.
\end{proo}

\subsection{Examples}  $n=1$: 

-- $\KK^{1} \qquad \KunAprim $
\\

-- fattened $ \KK^{1}$   $\qquad \Kungrossi$

-- Cone over fat-$\KK^{1}$ $=\KK^{2}\qquad \Kungrossicone$

Example $n=2$: 
\\

-- $\KK^{2}\qquad \KdeuxAprim$
\\

-- fattened  $\KK^{2}$ $\qquad \KdeuxAgrossi$
\\

-- Cone over fat-$\KK^{2}$ $=\KK^{3}\qquad \KtroisA$
\\

\section{Triangulation of the associahedron} From the construction  of $\KK^{n}$ out of $\KK^{n-1}$ performed in the preceding section, it is clear that one can triangulate  $\KK^n$ by induction. 

\subsection{Product of simplices} Let us recall that, if $\DD^n$ denotes the standard simplex, then, in the triangulation of $\DD^p\times \DD^q$, the simplices are indexed by the $(p,q)$-shuffles.

For instance the triangulation of $\DD^1\times \DD^1$ is:
$$\squareT$$
The triangulation of $\DD^2\times \DD^1$ (a prism) is made of three tetrahedrons:
$$\begin{array}{ccccc}
\textrm{shuffle} & \textrm{vertices} & & & \\
123 & (0,0), & (1,0), & (2,0), & (2,1)\\
132 & (0,0), & (1,0), & (1,1), & (2,1)\\
312 & (0,0), & (0,1), & (1,1), & (2,1)\\
\end{array}
$$
Here $(i,j)$ stands for the point of $\DD^2\times \DD^1$ which is the $i$th vertex of $\DD^2$ times the $j$th vertex of $\DD^1$.

\begin{thm}\label{simpl} The associahedron $\KK^n$ admits a triangulation by\\
 $(n+1)^{n-1}$ simplices, whose orientation is compatible with the Tamari order on the set of vertices. An $n$-simplex of this triangulation is a cone over an $(n-1)$-simplex determined by the following choices. Choose either

- an $(n-1)$-simplex of $\KK^{n-1}$, or

- in the fattened cell $\cc(a; p ,q)$ isomorphic to $\KK^{p }\times \KK^{q}$ choose a  $p$-simplex of $\KK^{p }$,  a $q$-simplex of $\KK^{q}$ and a $(p ,q)$-shuffle.
\end{thm}
\begin{proo} The proof of the second assertion follows from Theorem \ref{induction} and the fact that the triangulation of $\DD^{p }\times \DD^q$ is indexed by the $(p ,q)$-shuffles.

From this description of the triangulation we can count the number $d_n$ of top dimensional simplices of $\KK^n$  by induction. Let us suppose that $d_{p}= (p+1)^{p-1}$ for $p<n$, and recall that the number of $(p ,q)$-shuffles is the binomial coefficient
$\binom {p +q}{q}$. If $p $ is fixed, then $a$ can take the values $1,\ldots , p+1$. Hence there are $p+1 $ cells of the form $\KK^{p}\times \KK^{q}, p+q=n-1$. We get
\begin{eqnarray*}
d_n &=&  \sum_{p=0}^{n-1} (p+1) \binom{n-1}{p} d_{p}d_{n-p-1}\\
&=&  \sum_{p=0}^{n-1} \binom{n-1}{p} {(p+1)}^{p} (n-p)^{n-p-2}\\
&=& (n+1)^{n-1} .
\end{eqnarray*}
The last equality is a particular case of Abel's formula, cf.~\cite{Riordan},
$$x^{-1}(x+y+n)^n = \sum_{k=0}^{n}\binom{n}{k}(x+k)^{k-1}(y+n-k)^{n-k}$$
for $x=y=1$ and $k=n-p$. 
Observe that, at each step of the construction of $\KK^n$ out of $\KK^{n-1}$, the orientation of the simplices coincides with the orientation of the edges given by the Tamari poset order.
\end{proo}

\noindent {\bf Example:} triangulation of fat-$\KK^2$ giving rise to the triangulation (by tetrahedrons) of $\KK^3$:

$$\KtroisFtriang$$

\section{Parking functions}\label{parking} A \emph {parking function} is a sequence of integers $(i_1, \ldots , i_n)$ such that the associated ordered sequence $j_1\leq  \ldots \leq  j_n$ satisfies the following conditions: $1\leq j_k\leq k$ for any $k$. For instance there is only one parking function of length one : $(1)$, there are three parking functions of length two: $(1,2), (2,1), (1,1)$. There are sixteen parking functions of length three: the permutations of $(1,2,3), (1,1,3), (1,2,2), (1,1,2), (1,1,1)$ (remark that 6+3+3+3+1=16).
It is well-known, cf.~for instance \cite{NT}, that  there are $(n+1)^{n-1}$ parking functions of length $n$. We denote by $PF_{n}$ the set of parking functions of length $n$. We denote by $Sh(p,q)$ the set of permutations which are $(p,q)$-shuffles.


\begin{thm}\label{thm:shuffle} For any $n$ there is a bijection
$$\pi:  \bigcup_{{}^{p+q=n-1}_{p\geq 0, q\geq 0}}\{1,\ldots , p+1\}\times Sh(p,q)\times PF_{p}\times PF_{q}\longrightarrow PF_{n}$$
given by $\pi(a,\theta; f,g)= (a,\theta_{*}(f_{1},\ldots , f_{p}, p+1+g_{1}, \ldots , p+1+g_{q}))$.
\end{thm}
\begin{proo} Let $a\in \{1,\ldots , p+1\}, \theta\in Sh(p,q), f\in  PF_{p}$ and $g\in  PF_{q}$. 

Let us first show that the sequence 
$$x:= (a,\theta_{*}(f_{1},\ldots , f_{p}, p+1+g_{1}, \ldots , p+1+g_{q}))$$
 is a parking function. Let $(\phi_{1},\ldots , \phi_{p})$, resp.~$(\psi_{1},\ldots , \psi_{q})$, be the sequence $f$, resp.~$g$, put in increasing order. In the sequence $x$ put in increasing order we first find the sequence $f$ with the number $a$ inserted, then the sequence $\psi$. Since the sequence $\phi$ has $p$ elements and $1\leq a\leq p+1$ the expected inequality is true for $a$. It is also true for all the elements of $\phi$ since $\phi_{j}$ is either at the place $j$ or at the place $j+1$. The expected inequality is true for all the elements of the sequence $p+1+\psi$ since $p+1+\psi_{j}$ is at the place $p+1+j$. Hence $x$ is a parking function.
\smallskip

Let us now construct a map in the other direction. Let 
$${\underline a}=(a=a_1,a_{2},\ldots , a_{n})$$
 be a parking function, referred to as the original sequence. Let ${\underline x}=(x_{1},\ldots , x_{n})$ be the ordered sequence where $x_{j}=a$. Let $k$ be the smallest integer such that $k>j$ and $x_{k}= k$. Then we put $p+2=k$. It may happen that there is no such integer. In that case we put $p+2=n+1$, that is $p+1=n$ and so $q=0$. With these choices there exists $\ss \in S_{k-2}, \ss' \in S_{n-k+1}$ and $\theta \in Sh(k-2,n-k+1)$ such that $(a,\theta\circ (\ss\times \ss')(x_{1}, \cdots , x_{j-1}, x_{j}, \cdots , x_{k}, \cdots , x_{n}))$ is the original sequence ${\underline a}$.

Example: ${\underline a}= (3,6,1,7,2,1,3,6)$. Then we get ${\underline x}=(1,1,2,{\underline 3},3,6,6,7)$ (where $a$ has been underlined), $j=4$ and $k=6$, therefore $a=3, p=4, q=3$. The two parking functions are $(1,2,1,3)$ and $(1,2,1)$ and the $(4,3)$-shuffle is the permutation whose action on $u_{1}u_{2}u_{3}u_{4} v_{1}v_{2}v_{3}$ gives $v_{1}u_{1}v_{2}u_{2}u_{3}u_{4}v_{3}$.

Hence we have constructed a map from $PF_{n}$ to\\
 $\bigcup_{{}^{p+q=n-1}_{p\geq 0, q\geq 0}}\{1,\ldots , p+1\}\times Sh(p,q)\times PF_{p}\times PF_{q}$.
\smallskip

In order to show that it is the inverse of the previous map, it is sufficient to verify that our algorithm gives $k-2=p$ when we start with a parking function of the form $(a, \textrm{sh}(\textrm{pf}(1,\cdots ,p), \textrm{pf}(p+2,\cdots , p+1+q))$. First, in the ordered sequence of a parking function the first element is always $1$, hence $p+2$ is the smallest element in $\textrm{pf}(p+2,\cdots , p+1+q)$, or $q=0$ and $p=n+1$. Second, we know that $a\leq p+1$, so in the ordered sequence $a=x_{j}$ appears before $p+2=k$, hence $p+2 $ is at the place $p+2$, whence $k>j$ and we are done.
\end{proo}
\subsection{Remark} As a Corollary we get from Abel's formula (cf.~the proof of Theorem \ref{simpl}) the well-known result:
$$ \#PF_n =(n+1)^{n-1} .$$

\subsection{Examples}\ 
\medskip

\begin{tabular}{ c c c c   l }
$n$ & $a$ & $p$ & $q$ &  parking functions\\
\hline
1 & 1  & 0 & 0 & (1) \\
\hline
2 & 1 & 1 & 0 &  (1,1) \\
   & 1 & 0 & 1 &  (1,2) \\
   & 2 & 1 & 0 &  (2,1) \\
\hline
3 & 1 & 2 & 0 & (1,1,1)\quad (1,1,2)\quad (1,2,1)\\
   & 1 & 1 & 1 & (1,1,3)\quad (1,3,1) \\
   & 1 & 0 & 2 & (1,2,2)\quad (1,2,3)\quad (1,3,2) \\
   & 2 & 2 & 0 & (2,1,1)\quad (2,1,2)\quad (2,2,1)\\
   & 2 & 1 & 1 & (2,1,3)\quad (2,3,1)\\
   & 3 & 2 & 0 & (3,1,1)\quad (3,1,2)\quad (3,2,1)\\
 \hline
 \end{tabular}
 
\medskip
In the following statement we use Theorem \ref{simpl}.

\begin{thm}\label{thm:bij} Let $\ss$ be a simplex of $\KK^n$ determined either by

-- a simplex $\oo$ of $\KK^{n-1}$, or by

--  a triple $(a,p,q)$ (determining a face $\KK^p\times \KK^q$), where $1\leq a\leq p+1$, $p+q=n-1$, and a simplex $\aa$ of $\KK^p$, a simplex $\bb$ of $\KK^q$ and a $(p,q)$-shuffle $\theta$. 

The map $\Phi_n$, which assigns to $\ss$ the parking function

 $\Phi_{n}(\ss):= \big(1, 1+\Phi_{n-1}(\theta)\big)$ in the first case,

 $\Phi_{n}(\ss):= \big(a, \theta_{*}(\Phi_{p}(\aa),p+1+\Phi_{q}(\bb))\big)$ in the second one,

\noindent is a bijection from the set of $n$-simplices of $\KK^n$ to the set $PF_{n}$ of parking functions.
\end{thm}
\begin{proo} We work by induction on $n$. For $n=1$ there is no choice: $\Phi(1\textrm{-cell}) = (1)$. In the description of the triangulation of $\KK^n$ given in Theorem \ref{simpl} we have constructed a bijection between the set of $n$-simplices of $\KK^n$ and the set
$  \bigcup_{{}^{p+q=n-1}_{p\geq 0, q\geq 0}}\{1,\ldots , p+1\}\times Sh(p,q)\times PF_{p}\times PF_{q}$. By Theorem \ref{thm:shuffle} this set is in bijection with $PF_{n}$. It is immediate to check that the composite of the two bijections is the map $\Phi_{n}$ described in the statement of the Theorem.
\end{proo}

\subsection{Examples}

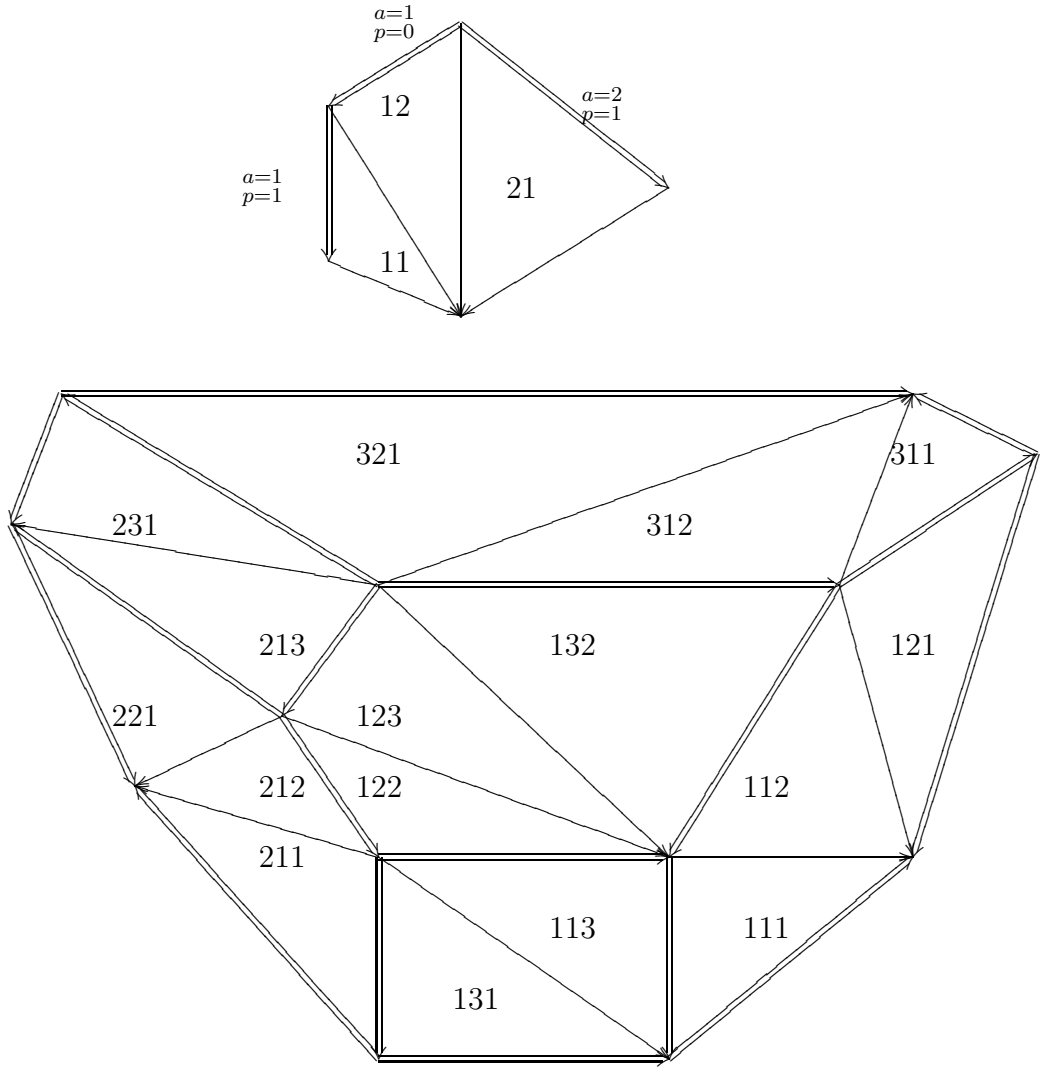
$$\KdeuxFindicbis$$

$$\KtroisFindic$$
\\

\subsection{Remark}\label{simplicial} Another way of stating Theorem \ref{simpl} and Theorem  \ref{thm:bij} is to say that we have constructed a simplicial set $K^n_{\cdot}$ such that $K^n_{0}= Y_{n}$,
$K^n_{n}/\{\textrm{degenerate elements}\}=PF_{n}$, and such that the geometric realization is 
$\vert K^n_{\cdot}\vert = \KK^n$ .

\subsection{Associahedron and cube}  It is well-known that for the cube $I^n$ the triangulation is indexed by the permutations, elements of the symmetric group $S_n$. In \cite{Lod2} we observed that the associahedron is contained in a certain cube. Locally around the North Pole (vertex with coordinates $(1,2,\ldots, n)$), the triangulation of the cube and the triangulation of the associahedron coincide. Our indexing of the simplices of the associahedron is such that the two indexings also coincide. 

\subsection{Other relationship between parking functions and associahedron} There are other links between parking functions and associahedron which are treated in \cite{PS} by Pitman and Stanley, in \cite{Po} by Postnikov and by Hivert (personal communication). It would be interesting to compare all of them.

\section{Triangulation of the permutohedron} In this section we briefly indicate how to simplicialize the permutohedron along the same line as the associahedron. So far we do not know of a nice combinatorial object playing the role of the parking functions.

\subsection{Permutohedron} Let us recall that the \emph{permutohedron} $\PP^{n-1}$ is the convex hull of the points $M(\ss)= (\ss(1), \ldots , \ss(n))\in \RRR^n$, where $\ss\in S_n$ is a permutation. 
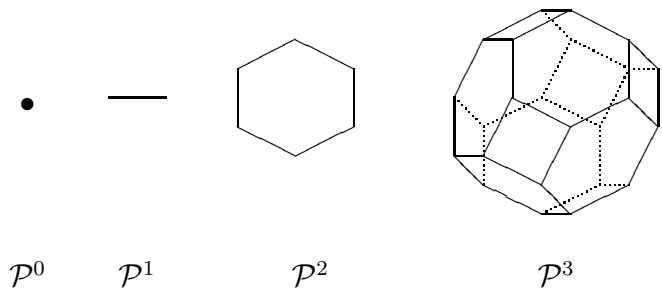
\begin{displaymath}
\begin{array}{cccc cccc}
&\Kzero &&  \KunA &&  \PdeuxA &&  \PtroisA \\
&&&&&&&\\
& \PP^0 && \PP^1 && \PP^2 && \PP^3\\
\end{array}
\end{displaymath}

The weak Bruhat order on $S_n$ is a partial order whose covering relations are in one to one correspondence with the edges of $\PP^{n-1}$. 

It is helpful to replace the permutations by \emph{planar binary trees with levels}. See for instance \cite{LR1} for a discussion of this framework. Under this remplacement the faces of the permutohedron $\PP^{n}$ are labelled by the planar leveled trees with $n+2$ leaves which have only two levels. Let us denote by $[p]$ the totally ordered set $\{0,\ldots , p\}$. A planar leveled tree with $n+2$ leaves and two levels is completely determined by the arity of the root, let us say $p+1$, and an ordered surjective map $f:[n+1] \to [p+1]$, $0\leq p\leq n-1$. We denote the corresponding face by $\gamma(p;f)$.\\

Example of faces of $\PP^2$:

\medskip

\begin{tabular}{rcl}
$\gamma(0; f)$ &=&$ \vcenter{\arbreAAC}$ where the image of $f$ is  (0,0,0,1),\\
$\gamma(0; f)$ &=& $\vcenter{\arbreACA}$ where the image of $f$ is  (0,0,1,1),\\
$\gamma(1;f)$ &=& $\vcenter{\arbreACC}$ where the image of $f$ is  (0,0,1,2).
\end{tabular}

\bigskip

The vertex $M(n+1, n, \ldots, 1)$, whose corresponding tree is the right comb, is called the South pole. For a given $p$ the face corresponding to the surjective map $f_0$, whose image is $(0, 1, \ldots, p, p, \ldots, p)$, contains the South pole.  We denote by $SM(n,p)$ the set of ordered surjective maps from $[n+1]$ to $[p+1]$ minus the map $f_0$. For instance $SM(4,1)$ has $10-1=9$ elements.

\subsection{Triangulation of $\PP^{n}$} We construct a simplicialization of  $\PP^{n}$ by induction as follows. For $n=0$, the space  $\PP^{0}$ is a point (0-simplex).  For $n=1$, since  $\PP^{1}$ is the interval (1-simplex). The permutohedron $\PP^{n}$ is the cone over the South pole with basis the union of the faces which do not contain the South pole, that is the faces whose ordered surjective map $f:[n+1] \to [p+1]$ is not the map $f_0$ given by $(0,1,\ldots , p+1,p+1,\ldots, p+1)$. Therefore an $n$-simplex of the triangulation of $\PP^{n}$ is the cone over the South pole $S$ with basis an $(n-1)$-simplex of the triangulation of the union of faces not containing $S$. Such an  $(n-1)$-simplex is completely determined by the following choices\\

- a face $\gamma(p ; f)$ not containing $S$ (i.e.~$f\neq f_0$),

- a top dimensional simplex in $\PP^{p}$,

- a top dimensional simplex in $\PP^{n-p-1}$,

- a shuffle in $Sh(p, n-p-1)$. \\

From this choice, it is clear by induction that the orientation of the simplices are compatible with the orientation of the edges induced by the weak Bruhat order.

Example: triangulation of $\PP^2$:
$$\Pdeuxoriente$$

In conclusion we have proved the following result.

\begin{thm}\label{zp} The set $ZP_n$ of top dimensional simplices ($n$-simplices) of the permutohedron $\PP^{n}$ satisfies the following recursive formula
$$ ZP_n = \bigcup_{p=0}^{n-1} SM(n,p) \times Sh(p, n-p-1)\times ZP_p \times ZP_{n-p-1},$$
where $ SM(n,p)= \binom{n+1}{p+1}-1$, and $ Sh(p,n-p-1)= \binom{n-1}{p}$.
\end{thm}
$\square$

The number of top simplices is as follows in low dimension:\\

\begin{tabular}{| c || c | c |  c | c | c | c |c |c |c |}
\hline
$n$ & 0& 1 & 2 & 3 & 4 & 5 & 6 & 7 & 8\\
\hline
$\#ZP_n$ & 1 & 1 & 4 & 34 & 488 & 10512 & 316224& 12649104 & 649094752\\
\hline 
\end{tabular}

\bigskip

\subsection{Permutohedron analogue of parking functions} It would be interesting to find a sequence of combinatorial objects analogue of the parking functions, that is satisfying the inductive relation of Theorem \ref{zp}.


\end{document}